# Finite dimensional factor algebras of $\mathbb{F}_2[X_1, \ldots, X_n]$ and their fixed point subalgebras

Miroslav Kureš
Institute of Mathematics
Brno University of Technology
Brno, Czech Republic
kures@fme.vutbr.cz

*Abstract*—Fixed point subalgebras of finite dimensional factor algebras of algebras of polynomials in *n* indeterminates over the finite field $\mathbb{F}_2$ (with respect to all $\mathbb{F}_2$-algebra automorphisms) are fully described.

## I. Introduction

In [6], we consider local commutative $\mathbb{R}$-algebra $A$ with identity, the nilpotent ideal $\mathfrak{n}_A$ of which has a finite dimension as a vector space and $A / \mathfrak{n}_A = \mathbb{R}$ and study its subalgebra $SA$ of fixed elements, $SA = \{a \in A;\ \varphi(a) = a$ for all $\varphi \in \mathrm{Aut}_\mathbb{R} A\}$, where $\mathrm{Aut}_\mathbb{R} A$ is the group of $\mathbb{R}$-automorphisms of the algebra $A$. This research is motivated by differential geometry, where algebras in question are usually called *Weil algebras* and, in particular, the bijection between all natural operators lifting vector fields from $m$-dimensional manifolds to bundles of Weil contact elements and the subalgebra of fixed points $SA$ of a Weil algebra $A$ was determined (in [4]). Although in the known geometrically motivated examples is usually $SA = \mathbb{R}$ (such $SA$ is called *trivial*), there are some algebras for which $SA$ is a proper superset of $\mathbb{R}$ and they call attention to the geometry of corresponding bundles. Thus, the fundamental problem is a classification of algebras having $SA$ nontrivial. See [4], [5] for related geometric questions and the survey paper [8] for known results up to now, especially for a number of claims concerning the form of subalgebras of fixed points of various Weil algebras.

In this paper, we simply replace $\mathbb{R}$ by $\mathbb{F}_2$ and study quite analogous questions. We come to a different situation by this: mainly, factor rings are <u>finite</u> rings (see [1]) and there is the whole theory about this topic. It is known the *ring automorphism problem* lying in a decision if a finite ring has a non-identical automorphism or not. Results about fixed point subalgebras are also qualitatively totally different from the real case and they can have interesting applications in the coding theory and cryptography.

## II. Polynomials over $\mathbb{F}_2$

Polynomials in $n$ indeterminates over $\mathbb{F}_2$ are the maps of the type

$$a: \mathbb{N}_0^n \to \mathbb{F}_2,$$

i.e. a multiindex maps onto an element of $\mathbb{F}_2$; the support of the map must be by definition finite. We define the addition of polynomials and the multiplication of polynomials by the usual way and use also the standard denotation for them. However, we consider only polynomials in this paper and not their evaluations (polynomial maps); the fatal inaccuracy of such a confusion is explained e.g. in [7]. With the mentioned operations, polynomials over $\mathbb{F}_2$ form the ring denoted by $\mathbb{F}_2[X_1, \ldots, X_n]$ or shortly by $\mathbb{F}_2[\mathbf{X}]$.

### II.1. Ideals in $\mathbb{F}_2[\mathbf{X}]$

The (unique) maximal ideal of $\mathbb{F}_2[\mathbf{X}]$ is

$$\mathfrak{m} = (X_1, \ldots, X_n).$$

Powers of $\mathfrak{m}$ represent notable class of ideals.

We mention another important ideal. For finite fields $\mathbb{F}_q$, the field ideal in $\mathbb{F}_q[X_1, \ldots, X_n]$ is defined as

$$\mathfrak{f} = (X_1^q - X_1, \ldots, X_n^q - X_n).$$

Thus, we have

$$\mathfrak{f} = (X_1^2 + X_1, \ldots, X_n^2 + X_n).$$

for $q = 2$.

## III. Algebras $(\mathbb{D}_2)_n^r$

In this section, we will study factor rings

$$(\mathbb{D}_2)_n^r = \mathbb{F}_2[X_1, \ldots, X_n] / \mathfrak{m}^{r+1},$$

where $r \in \mathbb{N}$.

### III.1. Dual numbers over $\mathbb{F}_2$

As $\mathbb{D} = \mathbb{R}[X] / (X^2)$ is usually called the *algebra of dual numbers* (which is definable promptly by $\mathbb{D} = \{a_0 + a_1 X;\ a_0, a_1 \in \mathbb{R},\ X^2 = 0\}$), we obtain for $r = n = 1$

$$\mathbb{D}_2 = (\mathbb{D}_2)_1^1 = \mathbb{F}_2[X] / (X^2)$$

the *algebra of dual numbers over* $\mathbb{F}_2$. Elements of $\mathbb{D}_2$ are expressible in the form

$$a_0 + a_1 X;\ a_0, a_1 \in \mathbb{F}_2,\ X^2 = 0.$$

---

Published results were acquired using the subsidization of the GA ČR, grant No. 201/09/0981.

We observe that $\mathbb{D}_2$ has the following additive and multiplicative tables:

| + | 0 | 1 | $X$ | $1+X$ |
|---|---|---|---|---|
| 0 | 0 | 1 | $X$ | $1+X$ |
| 1 | 1 | 0 | $1+X$ | $X$ |
| $X$ | $X$ | $1+X$ | 0 | 1 |
| $1+X$ | $1+X$ | $X$ | 1 | 0 |

| × | 0 | 1 | $X$ | $1+X$ |
|---|---|---|---|---|
| 0 | 0 | 0 | 0 | 0 |
| 1 | 0 | 1 | $X$ | $1+X$ |
| $X$ | 0 | $X$ | 0 | $X$ |
| $1+X$ | 0 | $1+X$ | $X$ | 1 |

As to classification of B. Fine, [2], this finite ring can be expressed as

$(a, b; 2a = 2b = 0, a^2 = 0, b^2 = b, ab = a, ba = a)$.

(the case denoted by 'G' in [2]).

Furthermore, in [3] is presented that dual numbers over $\mathbb{F}_2$, complex numbers over $\mathbb{F}_2$ and paracomplex numbers over $\mathbb{F}_2$ are isomorphic rings. We find easily the following result. (Analogously to the real case, by the *subalgebra of fixed elements SA* of an $\mathbb{F}_2$-algebra $A$ we mean the subalgebra of all elements $a$ satisfying $\varphi(a) = a$ for all $\mathbb{F}_2$-automorphisms $\varphi$ of $A$ and if $SA = \mathbb{F}_2$, we call $SA$ trivial.)

PROPOSITION 1. *The group of all $\mathbb{F}_2$-automorphisms of $\mathbb{D}_2$ is trivial. It follows that for $A = \mathbb{D}_2$ is SA nontrivial as $SA = A$.*

PROOF. Let $\varphi$ be an $\mathbb{F}_2$-automorphism of $\mathbb{D}_2$. As $\varphi(1) = 1$, $\varphi$ is fully determined by a specification of $\varphi(X)$. In general,

$\varphi(X) = b_0 + b_1 X;\; b_0, b_1 \in \mathbb{F}_2$.

However, we have

$0 = \varphi(0) = \varphi(X^2) = \varphi(X)\,\varphi(X) = b_0^2 + b_1^2 X^2 = b_0^2$;

thus, $b_0 = 0$, then, necessarily, $b_1 = 1$ for $\varphi$ be a bijection. So, the group of all $\mathbb{F}_2$-automorphisms of $\mathbb{D}_2$ contains only one element: the identical automorphism. Then the rest of the claim becomes evident. □

*III.2. The case $r>1$, $n=1$*

Elements of the algebra $(\mathbb{D}_2)^r = (\mathbb{D}_2)_1^r = \mathbb{F}_2[X] / \mathfrak{m}^{r+1} = \mathbb{F}_2[X] / (X^{r+1})$ have a form

$a_0 + a_1 X + a_2 X^2 + \ldots + a_r X^r;\; a_0, a_1, a_2, \ldots, a_r \in \mathbb{F}_2,\; X^{r+1} = 0$.

We start with the following lemma.

LEMMA 1. *Every endomorphism $\varphi: (\mathbb{D}_2)^r \to (\mathbb{D}_2)^r$ determined by*

$\varphi(1) = 1$
$\varphi(X) = X + b_2 X^2 + \ldots + b_r X^r;\; b_2, \ldots, b_r \in \mathbb{F}_2$

*belongs to the group of all $\mathbb{F}_2$-automorphisms of $(\mathbb{D}_2)^r$.*

PROOF. It suffices to describe $\varphi^{-1}$: we have

$Y = \varphi(X) = X + b_2 X^2 + \ldots + b_r X$
$Y^2 = X^2 +$ terms of degree $> 2$
...
$Y^{r-1} = X^{r-1} +$ terms of degree $> r-1$
$Y^r = X^r$

The last equation provides $X^r$ by $Y$'s, the last but one provides (after the substitution) $X^{r-1}$ and so on. □

We recall that $\mathfrak{n}_A$ denotes the ideal of nilpotent elements of $A$ (nilradical of $A$). If an element $a \in A$ has the property $au = 0$ for all $u \in \mathfrak{n}_A$, we call $a$ the *socle element* of $A$. It is easy to find that all socle elements constitute an ideal; this ideal is called the *socle* of $A$ and denoted by $\mathrm{soc}(A)$. Now, we can formulate the main result about automorphisms of $(\mathbb{D}_2)^r$ including also the case $r = 1$.

PROPOSITION 2. *For $r \in \mathbb{N}$, let $A = (\mathbb{D}_2)^r$. Every automorphism $\varphi: A \to A$ has a form*

$\varphi(1) = 1$
$\varphi(X) = X + b_2 X^2 + \ldots + b_r X^r;\; b_2, \ldots, b_r \in \mathbb{F}_2$.

*It follows SA is always nontrivial, in particular $\mathrm{soc}(A) \subseteq SA$.*

PROOF. It is evident that the endomorphism

$\varphi(1) = 1$
$\varphi(X) = b_2 X^2 + \ldots + b_r X^r;\; b_2, \ldots, b_r \in \mathbb{F}_2$

does not represent an automorphism. Further,

$\mathrm{soc}((\mathbb{D}_2)^r) = \{aX^r;\; a \in \mathbb{F}_2\}$

and $\varphi(X^r) = X^r$ was demonstrated already in the previous lemma. □

*III.3. The case $r=1$, $n>1$*

Elements of the algebra $(\mathbb{D}_2)_n = (\mathbb{D}_2)_n^1 = \mathbb{F}_2[X_1, \ldots, X_n] / \mathfrak{m}^2 = \mathbb{F}_2[X_1, \ldots, X_n] / (X_1, \ldots, X_n)^2$ have a form

$a_0 + a_1 X_1 + a_2 X_2 + \ldots + a_n X_n;\; a_0, a_1, a_2, \ldots, a_n \in \mathbb{F}_2,\; X_i X_j = 0$ for all $i, j \in \{1, \ldots, n\}$.

PROPOSITION 3. *The group of all $\mathbb{F}_2$-automorphisms of $(\mathbb{D}_2)_n$ is isomorphic to the general linear group $\mathrm{GL}(n,2) = \mathrm{GL}(n, \mathbb{F}_2)$ of the order $n$ over $\mathbb{F}_2$.*

PROOF. A general form of endomorphisms of $(\mathbb{D}_2)_n$ is

$\varphi(1) = 1$
$\varphi(X_1) = b_{10} + b_{11} X_1 + b_{12} X_2 + \ldots + b_{1n} X_n$
$\varphi(X_2) = b_{20} + b_{21} X_1 + b_{22} X_2 + \ldots + b_{2n} X_n$
...
$\varphi(X_n) = b_{n0} + b_{n1} X_1 + b_{n2} X_2 + \ldots + b_{nn} X_n$.

However, we have

$0 = \varphi(0) = \varphi(X_1^2) = \varphi(X_1)\,\varphi(X_1) = b_{10}^2 + b_{11}^2 X_1^2 + \ldots + b_{1n}^2 X_n^2 = b_{10}^2$,

thus, $b_{10} = 0$, and analogously, $b_{20} = \ldots = b_{n0} = 0$. Now, the matrix $(b_{ij})$ must be invertible for $\varphi$ be a bijection. So, automorphisms of $(\mathbb{D}_2)_n$ corresponds exactly with the group $\mathrm{GL}(n,2)$. □

COROLLARY 1. *Let $i, j \in \{1,...,n\}$, $i \neq j$. Then $\varphi_{(i,j)} : (\mathbb{D}_2)_n \to (\mathbb{D}_2)_n$ given by*

$$\varphi_{(i,j)}(1) = 1$$
$$\varphi_{(i,j)}(X_i) = X_i + X_j$$
$$\varphi_{(i,j)}(X_k) = X_k \text{ for all } k \in \{1,...,n\}, k \neq i,$$

*belongs to the group of all $\mathbb{F}_2$-automorphisms of $(\mathbb{D}_2)_n$.*

PROOF. It is clear that $\varphi_{(i,j)}$ meets the form from the proof of the previous proposition. □

REMARK 1. We remark that the order of GL($n$, 2) is

$$\Pi_{i=0}^{n-1}(2^n - 2^i).$$

PROPOSITION 4. *For $n \in \mathbb{N}$, $n > 1$, let $A = (\mathbb{D}_2)_n$. Then the subalgebra SA of fixed points of A is always trivial.*

PROOF. First, we prove that the element

$$X_1 + X_2 + \ldots + X_n$$

is not fixed. For this, it suffices to take some automorphism $\varphi_{(i,j)}$, e.g. $\varphi_{(1,2)}$ sends $X_1 + X_2 + \ldots + X_n$ onto $X_1 + X_3 + \ldots + X_n$. Second, let $\{k_1, \ldots, k_h\}$ be a (non-empty) proper subset of $\{1, ..., n\}$, i.e. $h < n$. We prove that the element

$$X_{k_1} + X_{k_2} + \ldots + X_{k_h}$$

is not fixed, too. We take $i \in \{k_1, \ldots, k_h\}$ and $j \in \{1, ..., n\} - \{k_1,...,k_h\}$ and apply $\varphi_{(i,j)}$: it sends $X_{k_1} + X_{k_2} + \ldots + X_{k_h}$ onto $X_{k_1} + X_{k_2} + \ldots + X_{k_h} + X_j$. So, $SA = \mathbb{F}_2$.

### III.4. The case $r > 1$, $n > 1$

Elements of the algebra $(\mathbb{D}_2)_n^r = \mathbb{F}_2[X_1, \ldots, X_n] / \mathfrak{m}^{r+1} = \mathbb{F}_2[X_1, \ldots, X_n] / (X_1, \ldots, X_n)^{r+1}$ have a form

$a_0 +$
$a_1 X_1 + a_2 X_2 + \ldots + a_n X_n +$
$a_{11} X_1^2 + a_{12} X_1 X_2 + \ldots + a_{nn} X_n^2 +$
$\ldots +$
$a_{1\ldots1} X_1^r + a_{1\ldots12} X_1^{r-1} X_2 + \ldots + a_{n\ldots n} X_n^r,$
$a_0, a_1, \ldots, a_n, a_{11}, \ldots, a_{n\ldots n} \in \mathbb{F}_2.$

On basis of previous results we can find out nature of this general case now.

PROPOSITION 5. *For $r \in \mathbb{N}$, $n \in \mathbb{N}$, $n > 1$, let $A = (\mathbb{D}_2)_n^r$. Then the subalgebra SA of fixed points of A is always trivial.*

PROOF. Obviously, elements of GL($n$,2) represent automorphisms also for $(\mathbb{D}_2)_n^r$. Of course, not <u>all</u> automorphisms, however, these (linear) automorphisms suffice for our following considerations. In the proof, we use formally partial derivations $\partial/\partial X_j$ for an expressing whether elements of $A$ contain $X_j$ in some non-zero power or not.

Let $P \in A$ and let exist $i, j \in \{1,..., n\}$ such that $\partial P/\partial X_i \neq 0$ and $\partial P/\partial X_j = 0$. Analogously with the case $r = 1$, $n > 1$, we apply $\varphi_{(i,j)}$ for the demonstration that $P$ cannot be fixed.

So, let $Q \in A$ is not of such a type and let $\sigma$ be a permutation of $n$-tuple $(X_1, \ldots, X_n)$ for which $\sigma(Q) \neq Q$. As permutations of $(X_1, \ldots, X_n)$ are also elements of GL($n$,2), we find again that $Q$ cannot be fixed.

Therefore we take $R \in A$ such that $\partial P/\partial X_i \neq 0$ for all $i \in \{1,..., n\}$ and such that does not exist any permutation of $(X_1, \ldots, X_n)$ yielding a transformation of $R$. Nevertheless, a "symmetry" of $R$ will be again unbalanced by $\varphi_{(i,j)}$, e.g. $\varphi_{(1,2)}$. Hence we have an automorphism for which not even $R$ is fixed.

Thus, only zero degree elements of $A$ remain fixed with respect to all automorphisms: $SA$ is trivial. □

### III.5. Summary

The previous assertions provide the following summary theorem ($r, n \in \mathbb{N}$).

THEOREM. *The subalgebra SA of fixed points of $A = (\mathbb{D}_2)_n^r$ is nontrivial if and only if $n = 1$.*